\documentstyle[12pt,verbatim]{amsart}

\let\mathcal\cal
\newtheorem{theorem}{Theorem}[section]
\newtheorem{proposition}[theorem]{Proposition}
\newtheorem{lemma}[theorem]{Lemma}
\newtheorem{corollary}[theorem]{Corollary}

\newtheorem{sublemma}{Sublemma}[theorem]

\newtheorem{case}{Case}

\newtheorem{Claim}{Claim}

\newtheorem{qtheorem}{Theorem}

\theoremstyle{definition}

\newtheorem{definition}[theorem]{Definition}

\theoremstyle{remark}
\newtheorem{remark}{Remark}

\newcommand{\proof}{\begin{pf}}
\newcommand{\Proof}[1]{\begin{pf*}{Proof of #1}}
\newcommand{\eproof}{\end{pf}}
\newcommand{\Eproof}{\end{pf*}}
\newcommand{\sproof}[1]{\begin{pf*}{#1}}
\newcommand{\esproof}{\end{pf*}}

\newcommand{\arablabel}{
          \renewcommand{\labelenumi}{{\rm (\arabic{enumi})}}
          \renewcommand{\theenumi}{{\rm (\arabic{enumi})}}
                    }

\newcommand{\alabel}{
          \renewcommand{\labelenumi}{{\rm (\alph{enumi})}}
          \renewcommand{\theenumi}{{\rm (\alph{enumi})}}
                    }
\newcommand{\Alabel}{
          \renewcommand{\labelenumi}{{\rm (\Alph{enumi})}}
          \renewcommand{\theenumi}{{\rm (\Alph{enumi})}}
                    }
\newcommand{\rlabel}{
          \renewcommand{\labelenumi}{{\rm (\roman{enumi})}}
          \renewcommand{\theenumi}{{\rm (\roman{enumi})}}
                    }

\newcommand{\acal}{{\mathcal A}}
\newcommand{\bcal}{{\mathcal B}}
\newcommand{\ccal}{{\mathcal C}}
\newcommand{\dcal}{{\mathcal D}}
\newcommand{\ecal}{{\mathcal E}}
\newcommand{\fcal}{{\mathcal F}}
\newcommand{\gcal}{{\mathcal G}}
\newcommand{\hcal}{{\mathcal H}}
\newcommand{\ical}{{\mathcal I}}

\newcommand{\pcal}{{\mathcal P}}

\newcommand{\scal}{{\mathcal S}}

\newcommand{\wcal}{{\mathcal W}}

\newcommand{\setm}{\setminus}
\newcommand{\empt}{\emptyset}
\newcommand{\subs}{\subset}
\newcommand{\sups}{\supset}
\newcommand{\oo}{{{\omega}_1}}
\newcommand{\rest}{\lceil}
\newcommand{\dom}{\operatorname{dom}}
\newcommand{\ran}{\operatorname{ran}}
\def\<{\left\langle}
\def\>{\right\rangle}
\def\cf{\operatorname{cf}}

\def\OO{{\omega}}
\def\oo{\omega_1}
\def\br#1;#2;{\bigl[ {#1} \bigr]^ {#2} }
\def\bc#1;#2;{\bigl( {#1} \bigr)^ {#2} }

\def\ooseq#1;#2;{\< {#1}_{#2}:{#2}<\oo\>}
\def\ooset#1;#2;{\{ {#1}_{#2}:{#2}<\oo\}}
\def\seq#1;#2;#3;{\< {#1}_{#2}:{#2}<#3\>}
\def\set#1;#2;#3;{\{ {#1}_{#2}:{#2}<#3\}}
\def\oseq#1;#2;{\< {#1}_{#2}:{#2}<\OO\>}
\def\oset#1;#2;{\{ {#1}_{#2}:{#2}<\OO\}}
\def\oosequ#1;#2;{\< {#1}^{#2}:{#2}<\oo\>}
\def\oosetu#1;#2;{\{ {#1}^{#2}:{#2}<\oo\}}
\def\sequ#1;#2;#3;{\< {#1}^{#2}:{#2}<#3\>}
\def\setu#1;#2;#3;{\{ {#1}^{#2}:{#2}<#3\}}
\def\osequ#1;#2;{\< {#1}^{#2}:{#2}<\OO\>}
\def\osetu#1;#2;{\{ {#1}^{#2}:{#2}<\OO\}}

\def\force{\raisebox{1.5pt}{\mbox {$\scriptscriptstyle\|$}}
\mbox{$\!\mbox{---}$}}

\newcommand{\fn}{\operatorname{Fn}}
\newcommand{\HH}{\operatorname{H}}

\def\to{\longrightarrow}

\newcommand{\newcases}{\setcounter{case}{0}}

\newcommand{\cont}{2^{\omega}}
\def\fin#1;{\br #1;<{\omega};}

\newcommand{\gdot}{\dot g}
\newcommand{\fnn}{{\fn}({\omega}_3,{\omega}_1;{\omega}_1)}

\newcommand{\sxe}{s({\xi},{\eta})}
\newcommand{\axe}{{\alpha}({\xi},{\eta})}

\newcommand{\px}{p({\xi})}
\newcommand{\pe}{p({\eta})}
\newcommand{\qe}{q({\eta})}
\newcommand{\qep}{q({\eta}')}
\newcommand{\rxe}{r({\xi},{\eta})}
\newcommand{\hxe}{h({\xi},{\eta})}
\newcommand{\hxep}{h({\xi},{\eta}')}
\newcommand{\fxe}{f({\xi},{\eta})}
\newcommand{\fxpep}{f({\xi}',{\eta}')}
\newcommand{\fxep}{f({\xi},{\eta}')}
\newcommand{\fxepp}{f({\xi},{\eta}'')}
\newcommand{\dx}{{\delta}({\xi})}
\newcommand{\pdx}{p^{\Delta}({\xi})}
\newcommand{\qde}{q^{\Delta}({\eta})}
\newcommand{\rpxe}{r^+({\xi},{\eta})}
\newcommand\cupd{\stackrel{.}{\cup}}
\newcommand\gbar{\overline\gcal}
\newcommand\Hbar{\overline\hcal}
\newcommand\adot{{\dot A}}
\newcommand\bdot{{\dot B}}
\newcommand\boole[1]{{[[#1]]}}
\def\oot{{{\omega}_2}}
\def\con#1;{\mbox{\rm Con}$($#1$)$}

\newcommand{\sab}{{\sigma}_{{\alpha},{\beta} }}

\newcommand{\tp}{\operatorname{tp}}
\newcommand{\w}{\operatorname{w}}
\def\qstar#1;{q^*_{#1}}

\newcommand{\xaca}{X[\acal]}

\newcommand{\ooto}{\oo\times{\omega}}
\newcommand{\piz}{\pi_0}
\newcommand{\pio}{\pi_1}

\newcommand{\lez}{\le_0}
\def\afi#1;{\<A_{#1},f_{#1}\>}
\newcommand{\af}{\afi;}
\def\afdi#1;{\<A_{#1},f_{#1},d_{#1}\>}
\newcommand{\afd}{\afdi;}
\newcommand{\bge}{\<B,g,e\>}

\newcommand{\bg}{\<B,g\>}
\newcommand{\up}{U^p}
\newcommand{\uq}{U^q}

\newcommand{\xn}{X^{(n)}}

\newcommand{\kxua}{\operatorname{K_{x,U,{\alpha}}}}
\newcommand{\cax}{\operatorname{C}}

\newcommand{\wov}{\overline{w}}

\def\npr#1;#2;#3;#4;{\prod^{#1}_{#2,#3}#4}

\newcommand{\coo}{\ccal_{{\omega}_1}}

\def\noc#1;{P^{\text{L}}_{#1}} 

\newcommand{\ssff}{\sigma_f}

\author{I. Juh\'asz}
\address{Mathematical Institute of the Hungarian Academy of Sciences
\\Budapest,  Re\'altanoda u, 13-15, 1053, Hungary}
\email{juhasz@@math-inst.hu}
\author{Zs. Nagy}
\address{Mathematical Institute of the Hungarian Academy of Sciences
\\Budapest, Re\'altanoda u, 13-15, 1053, Hungary}
\email{zsiga@@math-inst.hu}
\author{L. Soukup}
\address{Mathematical Institute of the Hungarian Academy of Sciences
\\Budapest,  Re\'altanoda u, 13-15, 1053, Hungary}
\thanks{The third author  was supported by DFG (grant Ko 490/7-1) and 
Magyar Tudom\'any\'ert Alap\'{\i}tv\'any}
\email{soukup@@math-inst.hu}
\author{ Z. Szentmikl\'ossy}
\address{Department of Analysis, E{\"o}tv{\"o}s University
\\Budapest, M\'uzeum krt, 6--8, 1088, Hungary}
\email{szentmiklossy@@math-inst.hu}
\title{Intersection properties of open sets, II}
\thanks{The preparation of this paper was supported by the 
Hungarian National Foundation for Scientific Research grant no. l908.}


\begin{document}
\begin{abstract}
A topological space is called $P_2$ ($P_3$, $P_{<{\omega}}$) if and only if  it does not contain two (three, finitely many)
uncountable open sets with empty intersection.
We show that (i) there are  0-dimensional $P_{<{\omega}}$ spaces  of size
$2^{\omega} $, (ii) there are
 compact $P_{<{\omega}}$ spaces of size ${\omega_1}$, (iii) the existence of 
a $\Psi$-like examples for (ii) is independent of ZFC, (iv) it is consistent
that $2^{\omega} $ is as large as you wish but every first
countable (and so every compact) $P_2$ space has cardinality 
$\le{\omega_1}$.  
\end{abstract}
\maketitle 
\section{Introduction}
\label{sc:intr}

In this paper we continue the investigations started  
in \cite{HJ}. There   $P_2$ spaces, 
i.e. spaces having no two uncountable disjoint open subsets, 
were considered.
We solve some problems which were left open in that paper and
strengthen some of its results. 
First we introduce some strengthenings of notion $P_2$.

\begin{definition}
A topological space $X$ has property {\em $P_n$},
where $n$ is a natural number, if it is Hausdorff and given open sets $U_0$,
$U_1$, $\dots$, $U_{n-1}$  with $\cap\{U_i:i<n\}=\empt$ 
we have $|U_i|\le {\omega} $ for some $i<n$. We say $X$ is 
{\em $P_{<{\omega} }$} if and only if  it is $P_n$ for each $n<{\omega} $.
The space $X$ is called {\em strongly $P_2$} provided the intersection of 
two uncountable open sets is uncountable.
\end{definition}

Clearly $P_m$ spaces are $P_n$ for $n<m$ and 
strongly $P_2$ spaces have property $P_{<{\omega}}$.

In section \ref{sc:zfc} we construct   
$P_{<{\omega}}$-spaces: a 0-dimensional one of size $2^{\omega}$ and 
two locally compact, first countable examples 
of size ${\omega_1} $. 
One of the ZFC constructions of locally compact $P_{<{\omega}}$ spaces 
is due to S. Shelah \cite{Sh}
and it is included here with his kind permission.
  
In section \ref{sc:ma} we will  see why the construction of  a 
compact $P_3$ space  can be expected to be much harder than  that of a $P_2$
one: 
a $\Psi$-like example (see the definition \ref{df:psi} below), which 
worked in the $P_2$ case,  
can not be constructed for  the  $P_3$ case in ZFC.
On the other hand we   show (corollary \ref{tm:big_psi}) that it is  
consistent that
$2^{\omega} $ is as large as you wish and there is a 
$\Psi$-like $P_{<{\omega}}$-space
of size $2^{\omega} $.

\begin{definition}
\label{df:psi} {\bf (1)}
Given an almost disjoint family $\acal\subs \br {\omega};{\omega};$ 
we define the topological space
$\xaca$ as follows: its underlying set is 
${\omega}\cup\acal$, the elements of ${\omega} $ are isolated and the 
neighborhood base of $A\in\acal$ in $\xaca$ is 
$\bcal_A=\{A\cup\{A\}\setm n:n\in {\omega} \}$.
{\bf (2)} A topological space  is {\em $\Psi$-like} 
 if and only if it is homeomorphic to
$\xaca$ for some almost disjoint $\acal\subs \br {\omega};{\omega};$.
\end{definition}

It is straightforward that $\Psi$-like spaces are  locally compact,
first countable and 0-dimensional.

In section \ref{sc:part} we present the proof of a partition theorem
below the continuum which was used in \cite{HJ} to show the
it is consistent that $2^{\omega} $ is big but every first countable
(and so every compact) $P_2$ space has cardinality $\le{\omega_2}$. 

Finally, in section \ref{sc:cohen} we strengthen this result (and 
solve problem 10 of \cite{HJ}) by showing
that it is consistent that $2^{\omega} $ is as large as you wish but 
every first countable (and so every compact) $P_2$ space is of size
$\le {\omega_1} $.

We use  standard topological notation and terminology throughout, cf \cite{J}.

We would like to thank for the referee for his helpful
suggestions and comments.

\section{ZFC constructions of $P_{<{\omega}}$-spaces}
\label{sc:zfc}
\begin{proposition}
\label{tm:p2-0dim}
There are 0-dimensional $P_{<{\omega} }$ spaces of size
$2^{\omega}$.
\end{proposition}

\sproof{First proof of proposition \ref{tm:p2-0dim}}
We show that the space $X$ from \cite[theorem 3]{HJ} is in fact
$P_{<{\omega}}$. To start with, we recall its definition. Fix an
independent family $\fcal\subs\br {\omega};{\omega};$ of size
$2^{\omega} $. The underlying set of $X$ will be 
${\omega} \cup \fcal$, where the elements of ${\omega} $ are isolated.
The neighborhood base of $F\in\fcal$ will be 
$\bcal_F=\{V_F(\gcal):\gcal\in\br \fcal\setm\{F\};<{\omega};\}$, where
$V_F(\gcal)=\{F\}\cup F\setm \cup\gcal$. 
If $\bigcup\fcal={\omega}$ then $X$ will be 0-dimensional and $T_2$.

Let $U_0,\dots,U_{n-1}$ be uncountable open subsets of $X$. For each 
$i<n$ we can find $\{F^i_{\nu}:{\nu}<\oo\}\subs \fcal$ and 
$\{\gcal^i_{\nu}:{\nu}<\oo\}\subs \br \fcal;<{\omega};$ such that 
$F^i_{\nu}\notin \gcal^i_{\nu}$, the $F^i_{\nu}$ are all distinct and 
$V_{F^i_{\nu}}(\gcal^i_{\nu})\subs U_i$.

Define the set mapping $f:\oo\to\br \oo;<{\omega};$ by the stipulation
$$
f({\nu})=\{{\mu}:\exists j<n\ F^j_{\mu}\in\bigcup_{i<n}\gcal^i_{\nu}\}.
$$
Since the $F^j_{\mu}$ are different we have 
$|f({\nu})|\le n|\bigcup_{i<n}\gcal^i_{\nu}|$. By Hajnal's theorem 
on set mappings (\cite{H}) 
there is an uncountable $f$-free subset $I$ of $\oo$.
Let ${\nu}_0,\dots,{\nu}_{n-1}$ be different elements of $I$.
Then $\{F^j_{{\nu}_j}:j<n\}\cap\bigcup_{i<n}\gcal^i_{{\nu}_i}=\empt$,
so 
$$
\empt\ne\bigcap_{j<n}F^j_{{\nu}_j}\setm
\bigcap(\bigcup_{i<n}\gcal^i_{{\nu}_i})= \bigcap_{i<n}
V_{F^i_{{\nu}_i}}(\gcal^i_{{\nu}_i})\cap \omega
\subs\bigcap_{i<n}U_i,
$$ 
which completes the proof.
\esproof

\sproof{Second proof of proposition \ref{tm:p2-0dim}} 
Consider the space  $D\cup\{p\}$ where $|D|=\cont$, the elements of
$D$ are isolated, and  the co-countable subsets of $D$ forms the 
neighbourhood base of $p$. 

Let us recall the following theorem \cite[theorem 4.4.4]{vM}.
\begin{qtheorem} If  $X$ is a $P$-space (i.e. the intersection of countable 
many open sets is open) with $\w(X)\le\cont$
then $X$ can be embedded into ${\beta}{\omega}\setm {\omega}$.
\end{qtheorem}

According to this theorem  the space 
$D\cup\{p\}$ can be regarded  as a subspace of 
 ${\beta}{\omega}\setm {\omega}$.
Consider the subspace $Y={\omega}\cup D$ of ${\beta}{\omega}$. If
$U$ is an uncountable open subset 
of ${\beta}{\omega}$ with $|U\cap Y|>{\omega}$
then $|U\cap D|>{\omega}$, so $p\in \overline {U\cap D}$.
But ${\omega}$ is dense in ${\beta}{\omega}$, 
so $U\cap D\subs\overline{U\cap {\omega}}$.
Hence  $p\in\overline{U\cap {\omega}}$, i.e., $U\cap {\omega}\in p$.
(Remember that the elements of ${\beta}{\omega}\setm {\omega}$ are just the 
ultrafilters on ${\omega}$). But the intersection of finitely many elements of
$p$ is not empty, so $Y$ is $P_{<{\omega}}$.  
\esproof

\begin{proposition}\label{tm:po_oo}
There are  first countable,
locally compact, scattered $P_{<{\omega}}$ spaces  of size ${\omega}_1$.
\end{proposition}

We present two different examples. The first is due to S. Shelah
\cite{Sh}.

\sproof{First proof of proposition \ref{tm:po_oo}}
The underlying set of our space $X$ will be $\ooto$.
For $x\in\ooto$ write $x=\<\piz(x),\pio(x)\>$. For $A\subs\ooto$
and $i<2$ put $\pi_i(A)=\pi_i''A$. If ${\alpha}<\oo$ let
$Y_{\alpha}=\{{\alpha}\}\times {\omega}$ and 
$X_{{\alpha}}=
{\alpha}\times {\omega}$.

Let $P_0$ be the family of pairs $p=\af$ which satisfy (i)--(iv) below:
\begin{enumerate}\rlabel
\item $A\in\fin\ooto;$,
\item $f$ is a function, $f:A\times A\to 2$.
\end{enumerate}
To formulate (iii) and (iv) write $U(x)=\up(x)=\{y\in A:f(y,x)=1\}$
for $x\in A$. 
\begin{enumerate}\rlabel\addtocounter{enumi}{2}
\item $x\in U(x)$, $U(x)\setm \{x\}\subs X_{\piz(x)}$,
\item $x\in U(y)$ implies $U(x)\subs U(y)$.
\end{enumerate}

For $p=\af$ and  $q=\bg$ from $P_0$ let 
$$
\begin{array}{rcl}
p\lez q&\text{if and only if}&A\supset B,\\
&&f\supset g,\\
&&\forall x,y, z_0,\dots, z_{l-1}\in B\\
&&\text{\quad if\ }\uq(x)\cap \uq(y)=\bigcup_{j<l}\uq(z_j)\\
&&\text{\quad then\ } \up(x)\cap \up(y)=\bigcup_{j<l}\up(z_j).
\end{array}
$$

For $p\in P_0$ let 
$$\wcal(p)=\{\<x,s\>\in A^p\times\fin A^p;:s\subs X_{\piz(x)}
\}.
$$
For $w=\<x,s\>\in\wcal(p)$ let $b(w)=x$ and 
$\up(w)=\up(x)\setm\bigcup_{y\in s}\up(y)$.
Let 
$$
\dcal(p)=\{\<w_0,w_1\>\in\wcal(p)\times\wcal(p):
\piz(b(w_0))<\piz(b(w_1))\}.
$$

Now let $P$ be the family of triples $\afd$, where $\af\in P_0$,
$d$ is a function, $d:\dcal(\af)\to {\omega}$, such that 
\begin{enumerate}\arablabel
\item if $d(\<w_0,w_1\>)=d(\<w'_0,w'_1\>)$ then
$\piz(b(w_1))\ne \piz(b(w'_1))$,
\item if $d(\<w_0,w_1\>)=d(\<w'_0,w'_1\>)$, 
${\alpha}\in \piz(A)$ and
${\alpha}<\piz(b(w_0))\le \piz(b(w'_1))$ then 
$$
\up(w_0)\cap \up(w'_1)\cap Y_{\alpha}\ne\empt.$$
\end{enumerate}
Write $p=\afdi p;$ for $p\in P$.
If $p=\afd$ and $q=\bge$ are from $P$ take $p\le q$  if and only if 
$\af\lez \bg$ and $d\supset e$. Let $\pcal=\<P,\le\>$.
For ${\alpha}<\oo$ let $E_{\alpha}=\{p\in P:{\alpha}\in \piz(A^p)\}$ 
and $\ecal=\{E_{\alpha}:{\alpha}\in \oo\}$.

\begin{lemma}
$($ZFC$)$ There is an $\pcal$-generic filter over $\ecal$.
\end{lemma}

\proof
This follows from \cite{Sh2}, or Velleman's $\omega$-morass (a weak form
of Martin's axiom provable in ZFC, see \cite{V}) can be applied for
$\pcal$ and $\ecal$. We explain the second argument more explicitly.
By \cite[Theorem 3.4]{V} it is enough to prove that 
$\ecal$ is weakly ${\omega}$-indiscernible for $\pcal$ (see 
\cite[Definition 1.5]{V}). 

For $n<{\omega}$, ${\alpha}<\oo$ and an order-preserving 
function $f:n\to {\alpha}$ define $\ssff:P_n\to P_{\alpha}$ as follows.
For $x=\<\nu,m\>\in n\times {\omega}$ let $\ssff(x)=\<f(\nu),m\>$.
For $\<x,s\>\in (n\times {\omega})\times \fin n\times {\omega};$ 
let $\ssff(\<x,s\>)=\<\ssff(x),\ssff''s\>$.
For $p=\<A,f,d\>\in P$ take $\ssff(p)=\<\ssff''A,f',d'\>\in P$, where
$f'(\ssff(x),\ssff(y))=f(x,y)$ and $d'(\<\ssff(w),\ssff(w')\>)=d(\<w,w'\>)$.

We should check conditions \cite[Definition 1.5]{V}(1)--(5). All but
(5) are clear. So let $s<n<{\omega}$, $f=f(s,n)$ (i.e. $\dom(f)=n$,
$f(i)=i$ for $i<s$ and $f(i)=n+i-s$ for $i\ge s$) and $p\in P_n$.
We need to find a common extension of $p$ and $q=\ssff(p)$.

Let $B=A_p\cup A_q$ and
$g:B\times B\to 2$ be the function defined by the stipulation
$g^{-1}\{1\}=f_p^{-1}\{1\}\cup f_q^{-1}\{1\}$. 
Clearly $r=\<B,g\>$ is a common extension 
of $p^-=\afi p;$ and $q^-=\afi q;$
in $\pcal_0$. Now for each 
$\wov_0=\<w_0,w_1\>\in \dcal(p)\setm \dcal(q)$, 
$\wov_1=\<w'_0,w'_1\>\in \dcal(q)\setm \dcal(p)$ and ${\alpha}\in \piz(B)$
with $d_p(\wov_0)=d_q(\wov_1)$ and 
${\alpha}<\piz(b(w_0))\le \piz(b(w'_1))$
 pick 
a different element
$z_{\{\wov_0,\wov_1,{\alpha}\}}\in Y_{\alpha}\setm B$.
Extend $r$ to $r'=\<B',g'\>\in\pcal_0$ by adding the points 
$z_{\{\wov_0,\wov_1,{\alpha}\}}$ to $B$ in such a way that
for any $x\in B$ we have $z_{\{\wov_0,\wov_1,{\alpha}\}}\in 
U^{r'}(x)$ if and only if  $b(w_0)\in U^r(x)$ or
$b(w'_1) \in U^q(x)$. 
Finally, taking $d=d_p\cup d_q$  
we choose $d':\dcal(r')\to{\omega}$ such that 
$d'\sups d$, $\ran(d'\setm d)\cap \ran(d)=\empt$ and $d'\setm d$
is 1--1. Now  $r^*=\<B',g',d'\>$ is a common
extension of $p$ and $q=\ssff(p)$ in $P$.
\eproof

Let $\gcal$  be an $\ecal$-generic filter over $\pcal$. 
The set $A=\bigcup\{A_p:p\in\gcal\}$ is uncountable because
$\piz(A)=\oo$ by the genericity of $\gcal$. 
Take $U(x)=\bigcup\{\up(x):x\in A^p\}$ and
$$
\bcal_x=\{U(x)\setm \bigcup_{j<l} U(z_j):
l<{\omega},z_0,\dots z_{l-1}\in X_{<\piz(x)}\cap A\}.
$$
Let $X=\<A,{\tau}\>$, where $\bcal_x$ is the neighborhood base 
of $x$ in $X$.

It is straightforward that $X$ is scattered, 0-dimensional and locally compact.

Finally we show that  $X$ is strongly $P_2$.  
Let 
$$\wcal=\{\<x,s\>\in A\times\fin A;:s\subs X_{\piz(x)}\}.
$$
For $w=\<x,s\>\in \wcal$ put $U_w=U(x)\setm\bigcup\limits_{x'\in s}U(x')$.

Assume now that $V_0$ and $V_1$ are uncountable open subsets of $X$.
We can find a sequence
$\scal=\{\<w^{\alpha}_0,w^{\alpha}_1\>:{\alpha}<\oo\}\subs\br \wcal;2;$
such that 
\begin{enumerate}\rlabel
\item $U_{w^{\alpha}_i}\subs V_i$ for $i<2$ and ${\alpha}<\oo$,
\item ${\alpha}+1=\piz(b(w^{\alpha}_0))<\piz(b(w^{\alpha}_1))$.
\end{enumerate}
Let $d=\bigcup\{d_p:p\in\gcal\}$. Then  $\scal\subs\dom(d)$, so there is
an uncountable $I\subs \oo$ such that $d$ is constant on 
$\{\<w^{\alpha}_0,w^{\alpha}_1\>:{\alpha}\in I\}$.
Then, by the definition of $\pcal$, we have 
$U_{w^0_{\alpha}}\cap U_{w^1_{\beta}}\cap Y_{\alpha}\ne\empt$ for
${\alpha}<{\beta}\in I$.
So the intersection $V_0$ and $V_1$ is uncountable. 
\esproof

\sproof{Second proof of proposition \ref{tm:po_oo}}
The construction will be divided into two parts. 
First we introduce the notion of ${\alpha}$-good spaces and we show that 
$\oo$-good spaces are $P_{<{\omega}}$. Then 
in the second part we construct 
an $\oo$-good space in $ZFC$.
\begin{definition}
Let $X=\<{\nu}\times {\omega},{\tau}\>$ be a locally compact scattered topological space of height
${\nu}\le\oo$ such that the ${\alpha}^{\rm th}$ level of $X$
is just $X_{\alpha}=\{{\alpha}\}\times {\omega}$. Write 
$\xn={\nu}\times \{n\}$. We say that $X$ is {\em ${\nu}$-good} if and only if 
\begin{enumerate}\alabel
\item $\xn$ is a closed subspace of $X$ and 
the natural bijection between $\xn$ and ${\nu}$ is a homeomorphism,
\item for each limit ordinal ${\beta}<{\nu}$ and each ${\alpha}<{\beta}$ if
$x\in X_{\beta}$ and $U$ is an open neighborhood of $x$ then the set
$$
\kxua=\{n\in{\omega}: U\cap \xn\cap 
(\bigcup\limits_{{\alpha}\le {\gamma}<{\beta}}X_{\gamma})=\empt\}
$$
is finite.
\end{enumerate} 
\end{definition}

\begin{lemma}\label{lm:good}
An $\oo$-good space $X$ is strongly $P_2$.
\end{lemma}

\Proof{lemma \ref{lm:good}}
Write 
$X_{<{\alpha}}$ for $\bigcup_{{\nu}<{\alpha}}X_{\nu}$.
Let $U$ and $V$ be uncountable open subsets of $X$. 
Let 
$$
I=\{n\in {\omega}:\{{\alpha}:\<{\alpha},n\>\in U\}\text{ is stationary}\}.
$$
By (a) we have $|\bigcup\{\xn:n\in I\}\setm U|\le {\omega}$.
Since the ideal $NS_{\oo}$ is ${\sigma}$-complete, 
the set $A=\{{\alpha}<\oo:\exists n\in{\omega}\setm I\ \<{\alpha},n\>\in U\}$
is not stationary.
Hence there is ${\alpha}<\oo$
such that $U\cap X_{\alpha}\subs \{{\alpha}\}\times I$.
But
$X_{\alpha}$ is dense in $X\setm X_{<{\alpha}}$, so 
$U\cap X_{\alpha}$ is dense in $U\setm X_{<{\alpha}}$.
It follows that $I$ is infinite.

We show  
$(U\cap V)\setm X_{<{\alpha}}\ne \empt$  for each ${\alpha}<\oo$.

Indeed pick $y\in V\cap X_{\beta}$ where ${\alpha}<{\beta}$ is limit.
By (b) the set $K_{y,V,{\alpha}}$ is finite, so we can choose
$n\in I\setm K_{y,V,{\alpha}}$. Then $V\cap \xn\cap(X_{<{\beta}}\setm
X_{<{\alpha}})\ne \empt$. But $|\xn\setm U|\le {\omega}$, so 
$V\cap U\cap (X_{<{\beta}}\setm X_{<{\alpha}})\ne\empt$ provided that 
${\alpha}$ is large enough, 
which completes  the proof of the lemma.
\Eproof
So it is enough  to construct an $\oo$-good space.

\begin{lemma}\label{lm:induction}
If ${\nu}<\oo$ and $X=\<{\nu}\times {\omega},{\tau}\>$ is a
${\nu}$-good space then there is a  ${\nu}+1$-good space $Y$ such that 
$X$ is a subspace of $Y$.
\end{lemma}

\Proof{lemma \ref{lm:induction}}
First recall that $X$ is collectionwise 
normal because it is countable and regular.
During the construction we have to distinguish two cases. 

\newcases
\begin{case}
${\nu}={\mu}+1$.
\end{case}

Let $\{x_{n,i}:n,i<{\omega}\}$ be a 1--1 enumeration
of $\{{\mu}\}\times {\omega}$.  
The family $\{x_{n,i}:n,i<{\omega}\}$ is closed and discrete, so 
applying the collectionwise normality of $X$
we can find a closed discrete family $\{U_{n,i}:n,i<{\omega}\}$ of
 clopen subsets of $X$ with $x_{n,i}\in U_{n,i}$. 
Take $W_{n,i}=\{\<{\nu},n\>\}\cup\bigcup\{U_{n,j}:i\le j<{\omega}\}$.
Let the neighborhood base of $\<{\nu},n\>$ in $Y$ be 
$\wcal_n=\{W_{n,i}:i<{\omega}\}$. 

\begin{case}
${\nu}$ is limit.
\end{case}

Let $\{{\gamma}_n:n<{\omega}\}$ be a strictly increasing, cofinal
sequence in ${\nu}$.
By induction choose ordinals ${\delta}_n<{\nu}$ and distinct points
$\{x_{n,i}:n\le i<{\omega}\}\subs X$ such that 
\begin{enumerate}\rlabel
\item ${\gamma}_n< {\delta}_n$,
\item $x_{n,i}\in ({\delta}_n\setm {\gamma}_n)\times \{n\}$.
\end{enumerate}

Put $E_n=({\nu}\setm {\delta}_n)\times \{n\}$.
Then $\fcal=\{E_n:n\in {\omega}\}\cup\{x_{n,i}:n\le i\in {\omega}\}$
is a closed, discrete family of closed sets because for each ${\alpha}<{\nu}$
the set $\{F\in\fcal:F\cap X_{<{\alpha}}\ne\empt\}$ is finite. 
Applying the collectionwise normality of $X$ 
we can find a closed, discrete family
$\{U_n:n\in {\omega}\}\cup\{V_{n,i}:n\le i<{\omega}\}$ of clopen
subsets of $X$ such that $E_n\subs U_n$ and $x_{n,i}\in V_{n,i}$.
For each $x=\<{\alpha},m\>\in X$ fix   a countable  neighborhood base 
$\bcal_x$ of $x$ containing only clopen subsets of $X_{\le {\alpha}}$. 
The neighborhood base of $\<{\nu},n\>$ in $Y$ will be generated by the sets
$$
W_{n,j,y,B}=(\bigcup_{j\le i<{\omega}}V_{n,i})\cup(U_n\setm B)
$$
where $n\le j<{\omega}$, $y\in U_n$ and $B\in \bcal_y$.
It is easy to see that $Y$ satisfies the requirements
\Eproof
Now define, by induction on ${\nu}$,
${\nu}$-good spaces $X_{\nu}$ for ${\nu}\le \oo$ which extend each other. 
For successor ${\nu}$ we can apply  lemma \ref{lm:induction}.
For limit ${\nu}$ we can simply take the union.
\esproof

\section{$\Psi$-like spaces}
\label{sc:ma}
All examples of locally compact $P_2$ spaces  in \cite[theorem 3,4 and 8]{HJ}
are $\Psi$-like. In this section we will see why the examples 
of first countable, locally compact $P_{<\omega}$ spaces 
constructed in the previous section in ZFC are not $\Psi$-like. 

\begin{theorem}\label{tm:ma_nop}
If ${MA}_{{\omega}_1}(\text{${\sigma}$-linked})$ holds then there is no 
 $\Psi$-like $P_3$-space.
\end{theorem}

Theorem \ref{tm:ma_nop} follows from the following
combinatorial result.

\begin{theorem}\label{tm:ma_part}
$(MA_{\kappa}(\text{${\sigma}$-linked}))$ 
If $\{F_{\alpha}:{\alpha}<{\kappa}\}\subs 
\br {\omega};{\omega};$ is an almost disjoint family then 
${\omega}$ and ${\kappa}$ have partitions 
$(a_0,a_1,a_2)$ and $(\ical_0,\ical_1,\ical_2)$, respectively, with 
$|a_i|={\omega}$ and $|\ical_i|={\kappa}$ such that
\begin{equation}
\tag{$*$}\forall i\in 3\ \forall {\alpha}\in \ical_i
\ |F_{\alpha}\cap a_i|<{\omega}.
\end{equation}
\end{theorem}

\begin{remark}
As the referee pointed out, $MA(\text{$\sigma$-centered})$ does not imply 
this statement  because the strong Luzin property 
(see definition \ref{df:s-Luzin} below) of an almost disjoint family
is preserved by any $\sigma$-centered forcing.
As it was recalled in \cite{HJ} partitioning ${\omega} $ 
into two pieces is also not enough:
 a Luzin gap , i.e,  an almost disjoint family 
$\{B_{\beta}:{\beta}<{\omega}_1\}\subs \br {\omega};{\omega};$
such that for each $B\in \br {\omega};{\omega};$ either the set
$\{{\beta}<{\omega}_1:|B_{\beta}\cap B|<{\omega}\}$ or the set 
$\{{\beta}<{\omega}_1:|B_{\beta}\setm  B|<{\omega}\}$ is at most
countable, can be constructed in ZFC. 
\end{remark}

\Proof{theorem \ref{tm:ma_part}}
For $I\subs {\kappa}$ write $F_I=\bigcup_{{\alpha}\in I}F_{\alpha}$.
Define the poset $\pcal=\<P,\le\>$ as follows. Its underlying set $P$
consists of 7-tuples $\<m,A_i,I_i:i\in 3\>$
satisfying (i)--(iii) below:
\begin{enumerate}\rlabel
\item $m\in{\omega}$, $A_i\in\br {\omega};<{\omega};$,
$I_i\in \br {\kappa};<{\omega};$,
\item $(A_0,A_1,A_2)$ is a partition of $m$,
\item $F_{I_i}\cap F_{I_j}\subs m$ for $i\ne j$.
\end{enumerate}

Write $p=\<m^p,A^p_i,I^p_i:i\in 3\>$ and 
$I^p=\bigcup_{i\in 3}I^p_i$ for $p\in P$.
Given $p,q\in P$ we set $p\le q$ if and only if  
\begin{enumerate}\alabel
\item $m^q\le m^p$,
\item $A^q_i=A^p_i\cap m^q $  for $i\in 3$,
\item $I^q_i\subs I^p_i$ for $i\in 3$,
\item $F_{I^q_i}\cap ( A^p_i\setminus A^q_i)=\empt$ for $i\in 3$.
\end{enumerate}
Obviously $\le$ is a partial order on $P$.

For $m\in {\omega}$ and $i\in 3$ put 
$$D_{m,i}=\{p\in P:A^p_i\setm m\ne\empt\}
$$
and
$$
D_m=\{p\in P:m^p\ge m\}.
$$

\begin{lemma}\label{l:dmi}
$\forall m\in {\omega}$ $\forall i\in 3$
$\forall q$ $\exists p\in D_{m,i}$
$p\le q$ $\land $ $I^p=I^q$.
\end{lemma}

\proof
Let $q\in P$ with $A^q_i\subs m$. We can assume that $m\ge m^q$.
Pick $m^*\in ({\omega}\setm F_{I^q_i})\setm m$.
Choose a function $f:[m^q,m^*]\to 3$
such that $f(m^*)=i$ and 
$k\notin F_{I^q_{f(k)}}$ for $k\in [m^q,m^*)$. Such a function
exists by (iii).
Let $A_j=A^q_j\cup f^{-1}\{j\}$ for $j<3$ 
and $p=\<m^*+1,A_j,F^q_j:j\in 3\>$.
Obviously $p\in P$ and $m^*$ witnesses $p\in D_{m,i}$.
Moreover, $p\le q$. Indeed, conditions (a)--(c) are trivial
and (d) holds because $A^p_j\setm A^q_j=f^{-1}\{j\}$.
\eproof

\begin{lemma}\label{l:dm}
$\forall m\in {\omega}$ $\forall q$ $\exists p\in D_m$
$p\le q$ $\land $ $I^p=I^q$.
\end{lemma}

\proof
Straightforward by lemma \ref{l:dmi}.
\eproof

\begin{lemma}\label{l:i}
If ${\alpha}\notin I^q$ then for each $i\in 3$ there is a $p\le q$
with ${\alpha}\in I^p_i$.
\end{lemma}

\proof
We can assume that $i=0$. Fix $m\in {\omega}$ with 
$F_{\alpha}\cap F_{I^q}\subs m$. Using lemma \ref{l:dm}
pick $r\in D_m$ such that $r\le q$ and $I^q=I^r$.
Put $p=\<m^r,A^r_j:j\in 3, I^r_0\cup\{{\alpha}\},I^r_1,I^r_2\>$.
Then $p\in P$ because (iii) holds by
$F_{\alpha}\cap F_{I^p\setm \{{\alpha}\}}=F_{\alpha}\cap F_{I^q}\subs m
\subseteq m^r=m^p$.
Since $p\le q$ is trivial, the lemma is proved.
\eproof

For ${\alpha}\in {\kappa}$ put 
$$
E_{{\alpha},i}=\{p\in P:I^p_i\cap
[{\omega}{\alpha},{\omega}{\alpha}+{\omega})\ne\empt\}
$$
and 
$$
E_{\alpha}=\{p\in P:{\alpha}\in I^p\}.
$$

\begin{lemma}\label{l:eai}
Both $E_{{\alpha},i}$ and $E_{\alpha}$ are dense in $\pcal$.
\end{lemma}

\proof
Straightforward by lemma \ref{l:i}.
\eproof

\begin{lemma}\label{l:kprop}
$\pcal$ is ${\sigma}$-linked.
\end{lemma}

\proof
By ${MA}_{{\omega}_1}(\text{${\sigma}$-linked})$ we have
${\kappa}<\cont$ and so we can choose
a countable dense set  $D$   
in the product space $3^{\kappa}$.
For each $p=\<m,A,I_i:i\in 3\>\in \pcal$ pick $d_p\in D$ such that 
$d_p({\alpha})=i$ whenever $i\in 3$ and ${\alpha}\in I_i$.

Let   $p_0$ and $p_1$ be  elements of $P$,  
$p_j=\<m^j,A^j_i,I^j_i:i\in 3\>$. We show that if  
\begin{enumerate}\arablabel
\item $m^0=m^1$,
\item $A^0_i=A^1_i$ for $i\in 3$,
\item $d_{p_0}=d_{p_1}$,
\end{enumerate}
then $p_0$ and $p_1$ are compatible. Clearly this statement
yields that $\pcal$ is ${\sigma}$-linked.

Let $I^j=\bigcup\limits_{i\in3}I^j_i$ for $j\in 2$ and $I=I^0\cap I^1$.
Fix first $m'\ge m$ with
\begin{equation}
\tag{\dag}F_{I^0\setm I}\cap F_{I^1\setm I}\subs m'.
\end{equation}
If $k\in [m,m')$ and $j\in 2$ then, 
by (iii), there is at most one $i_j(k)\in 3$
with $k\in F_{I^j_{i_j(k)}}$.
Choose $g(k)\in 3\setm \{i_0(k), i_1(k)\}$.
Then $k\notin F_{I^0_{g(k)}\cup I^1_{g(k)}}$.

Define the partition $(A'_0,A'_1,A'_2)$ of $m'$
by the equations $A'_i=A_i\cup g^{-1}\{i\}$ for  $i<3$.

Put $p=\<m',A'_i,I^0_i\cup I^1_i:i\in 3\>$.

Check first $p\in P$. Conditions (i) and (ii) are trivial. 
Since $I^0_i\cap I=I^1_i\cap I$ by (3) we have 
$$
F_{I^0_i\cup I^1_i}\cap F_{I^0_j}=(F_{I^0_i\cap I}\cap F_{I^0_j})\cup
(F_{I^1_j\setm I}\cap F_{I^1_j})\cup (F_{I^0_i\setm I}\cap F_{I^1_j\setm I})
\subs m\cup m\cup m'=m'
$$ 
 by (iii) for $p_0$ and $p_1$ and by (\dag). 
Thus (iii) holds for $p$.
Finally we show that $p$
is a common extension of $p_0$ and $p_1$.
Conditions (a), (b) and (c) obviously hold.
But $A^p_i\setm A^{p_0}_i=A'_i\setm A_i=g^{-1}\{i\}$
and $g^{-1}\{i\}\cap F_{I^0_i\cup I^0_i}=\empt$ by
the choice of $g$, so (d) is also satisfied.
The lemma is proved.
\eproof

We are now ready to conclude the proof of the theorem.
Let 
$$
\dcal=\{D_m,D_{m,i}:m\in{\omega},i\in 3\}\cup
\{E_{\alpha},E_{{\alpha},i}:{\alpha}\in{\kappa},i\in 3\}.
$$
By lemma \ref{l:dmi}-\ref{l:eai}, $\dcal$ is a family
of dense subsets of $\pcal$.
Since $\pcal$ is ${\sigma}$-linked, 
${MA}_{{\kappa}}(\text{${\sigma}$-linked})$ 
implies that there is a
$\dcal$-generic filter $\gcal$ over $\pcal$.
Put $a_i=\bigcup\{A^p_i:p\in\gcal\}$ and 
$\ical_i=\bigcup\{I^p_i:p\in\gcal\}$.
Then (ii) and (b) imply $a_i\cap a_j=\empt$ for $i\ne j$.
Since $D_m\cap\gcal\ne\empt$, we have 
$a_0\cup a_1\cup a_2\supset m$. Moreover $D_{m,i}\cap \gcal\ne\empt$
implies $a_i\not\subset m$. Thus $(a_0,a_1,a_2)$
is a partition of ${\omega}$ into infinite pieces.
Similar arguments show that $(\ical_0,\ical_1,\ical_2)$
is a partition of ${\kappa}$ into subsets of size ${\kappa}$.

Finally let $i\in 3$ and ${\alpha}\in \ical_i$.
Pick $q\in \gcal $ with ${\alpha}\in I^q_i$.
Then $F_{\alpha}\cap a_i\subs F_{I^q_i}\cap a_i=
\bigcup\limits_{p\in\gcal,p\le q }(F_{I^q_i}\cap A^p_i)\subs A^q\subs m^q$
by (d).
\Eproof

By the previous theorem the existence of  a $\Psi$-like $P_3$ space
of size $\oo$ is not provable in ZFC. On the other hand, we will see
that it is also consistent that  
$2^{\omega} $ is arbitrarily large  and there is a 
$\Psi$-like $P_{<{\omega}}$-space
of size $2^{\omega} $.

We start with a definition.
\begin{definition}\label{df:s-Luzin}
Let $\acal \subs\br {\omega}; {\omega} ; $ be an almost disjoint
family. We say that $\acal$ is a {\em strong Luzin gap } 
 if and only if given finitely many uncountable subsets $\acal_0,\dots,\acal_{n-1}$
of $\acal$ we have that 
$\bigcap\limits_{i<n}\left(\bigcup\acal_i\right)$ is infinite.
\end{definition}

It is easy to see that $\acal$ is a strong Luzin gap if and only if  
the corresponding $\Psi$-like space $\xaca$ (see definition 
\ref{df:psi}) is $P_{<{\omega}}$.    
 
In  \cite[theorem 6]{HJ} a Luzin gap was obtained by a c.c.c forcing 
due to Hechler \cite{He}.

We show that this almost disjoint family is in fact strongly Luzin.
First we recall some notations and definitions from \cite{HJ}.
Let ${\kappa}>{\omega}$ be a fixed regular cardinal.
Write $\dcal=\br {\kappa};<{\omega};\times {\kappa}$.
Let
$$
P=\{p:p\text{ is a function},\dom(p)\in\dcal,\ran(p)\subs 2\}
$$
If $p,p'\in P$, $\dom(p)=A\times n$, $\dom(p')=A'\times n'$   put  
$p\le p'$ if and only if  $p\supset p'$ and for each $k\in n\setm n'$ we have
$|\{{\alpha}:p({\alpha},k)=1\}|=1$. 

If $\gcal$ is $P$-generic then let
$A_{\alpha}=\{k:\exists p\in \gcal\ p({\alpha},k)=1\}$.
Take $\acal[\gcal]=\{A_{\alpha}:{\alpha}<{\kappa} \}$.


\begin{lemma}\label{lm:strong_luzin}
 $\acal[\gcal]$ is a strong Luzin gap in $V^P$.
\end{lemma}

\proof
Assume on the contrary, that $n$ and $m$ are natural numbers,
$p\force$ ``{\em $\{{\dot\nu}_i:i<n\}$  are injective functions from
${\omega_1}$ to ${\kappa}$ and }
\begin{equation}\tag{$*$}
\left(\bigcap_{i<n}\bigcup_{{\alpha}<\oo}A_{{\dot\nu}_i({\alpha})}
\right)\subs m.\text{''}
\end{equation}

For each ${\alpha}<\oo$ pick $p_{\alpha}\le_P p$ and
${\nu}_{\alpha,i}\in {\omega_1} $ for $i<n$ such that 
\begin{equation}\notag
p_{\alpha}\force\text{``}{\dot\nu}_i({\alpha})={\hat\nu}_{\alpha,i}
 \text{\ for each $i<n$.''}
\end{equation} 

Without loss of generality the ${\nu}_{\alpha,i} $ 
are pairwise different. Let $\dom(p_{\alpha})=D_{\alpha}\times k_{\alpha}$. 
We can assume that $k_{\alpha}$ are independent of ${\alpha}$, 
$k_{\alpha}=k\ge m$, $\{D_{\alpha}:{\alpha}<{\omega_1} \}$ forms a
$\Delta$-system with kernel $D$, $|D_{\alpha}|=|D_{\beta}|$,
$\{{\nu}_{{\alpha},i}:i<n\}\subs D_{\alpha}\setm D$ and denoting by
$\sab$ the unique order preserving bijection between $D_{\alpha} $
and $D_{\beta} $ we have that $\sab$ gives an isomorphism between
$p_{\alpha} $ and $p_{\beta} $ and 
$\sab({\nu}_{\alpha,i} )={\nu}_{\beta,i} $  for $i<n$. 

Let 
${\alpha}_0<\dots {\alpha}_{n-1}<\oo$.
Define the condition $q\in P$ as follows:
\begin{enumerate}
\item $\dom(q)=(\bigcup_{i<n}D_{{\alpha}_i})\times k+1$,
\item $q\supset \bigcup\limits_{i<n}p_{\alpha_i}$,
\item $q({{\nu}_{\alpha_i,i}},k)=1$ for $i<n$,
\item $q({\xi},k)=0$ provided 
${\xi}\in \bigcup\limits_{i,n}D_{{\alpha}_i}\setm 
\{{\nu}_{{\alpha}_0,0},\dots,{\nu}_{{\alpha}_{n-1},n-1}\}$.
\end{enumerate} 

Then $q$ is a condition which extends  $p_{\alpha_i} $ for $i<n$, 
but $q\force k\in\bigcap_{i<n}A_{{\nu}_{\alpha_i,i} }$ which contradicts
$(*)$ because $k\ge m$.  
This concludes the proof of the lemma.
\eproof

This lemma yields the following corollary.
\begin{corollary}\label{tm:big_psi} If ZF is consistent then so is ZFC +
 $2^{\omega}={\kappa}$ is as large as you wish +  there is a  
$\Psi$-like  $P_{<{\omega}}$ space of size $2^{\omega}$.
\end{corollary}

\section{A positive partition theorem below $2^{\omega}$} 
\label{sc:part}
In this section we present the proof of a theorem of Szentmikl{\'o}ssy 
which was announced and applied in \cite{HJ}.

Given $A,B\subs {\kappa}$ we denote by $[A;B]$ the set
$\{\{{\alpha},{\beta}\}:{\alpha}\in A,{\beta}\in B,{\alpha}<{\beta}\}$.

\begin{definition}
Given cardinals ${\kappa}$, ${\lambda}$ and ${\mu}$ the partition
relation ${\kappa}\to \big(({\lambda};{\lambda})\big)^2_{\mu}$
holds if and only if  
$\forall f:\br {\kappa};2;\to {\mu}$
$\exists A,B\in\br {\kappa};{\lambda};$ $\exists {\xi}<{\mu}$
$\big( \tp(A)=\tp(B)={\lambda}$ $\land$ $\sup A\le\sup B$ $\land$   
$f''[A;B]=\{{\xi}\}\big)$.
\end{definition}

\begin{theorem}\label{tm:zoli}
If ZF is consistent then so is ZFC $+$ $2^{\omega}={\kappa}\ge{\omega}_3$ $+$
${\omega}_3\to (({\omega}_1;{\omega}_1))^2_{\omega_1}$.
\end{theorem}

The proof is based on the following lemma.
To simplify its formulation we write $P=\fnn$.
For $f$, $g$, $h\in P$ we write $f=g\cupd h$ to means 
that $f=g\cup h$ and $g\cap h=\empt$.
\begin{lemma}\label{lm:zoli}
Assume GCH. If $f:[{\omega}_3]^2\to P$ then there is some 
$C\in[{\omega}_3]^{\omega_2}$ and there are elements 
$\{p,\px,\qe,\rxe:{\xi},{\eta}\in C,{\xi}<{\eta}\}\subs P$ 
with the following properties:
\begin{enumerate}\rlabel
\item $\forall{\xi}<{\eta}\in C$ $\fxe=p\cupd \px\cupd\qe\cupd\rxe$,
\item $\forall {\xi}\in C$ $\{\fxe:{\eta}\in C\setm {\xi}+1\}$
forms a $\Delta$-system with kernel $p\cupd \px$,
\item $\forall {\eta}\in C$ $\{\fxe:{\xi}\in C\cap {\eta}\}$
forms a $\Delta$-system with kernel $p\cupd \qe$,
\item the $\px$ have disjoint supports,
\item the $\qe$ have disjoint supports.
\end{enumerate}
\end{lemma}

\Proof{lemma}
Fix a large enough regular cardinal ${\tau}$ and let $N$ be a countable
elementary submodel of $\hcal({\tau})=\<{\HH}({\tau}),\in,<\>$ with
$f\in N$, where ${\HH}({\tau})$ is the family of the sets whose transitive 
closure has cardinality $<{\tau}$.

Let $\gcal=\{g\in N:g$ is a function from $\br {\omega}_3;<{\omega};$
to $\oo$ $\}$ and $\hcal=\{h\in N:h $ maps $\br {\omega}_3;2;$ to $\oo\}$.

\begin{sublemma}
There is  $D\subs {{\omega}_3}$ of order type $\oot+1$ which is 
end-homogeneous for all $g\in\gcal$ and homogeneous for all $h\in \hcal$.
\end{sublemma}

\Proof{sublemma}
By $2^{\oo}=\oot$ we can choose an increasing sequence 
$\<N_{\alpha}:{\alpha}\le\oot\>$ of elementary submodels of $\hcal({\tau})$
such that $N\in N_0$, $|N_{\alpha}|=\oot$ and 
$\br N_{\alpha};\oo;\subs N_{\alpha}$.

Since $\gcal,\hcal\in N_0$, they have enumerations 
$\gbar=\<g_n:n\in {\omega}\>$ and $\Hbar=\<h_n:n\in {\omega}\>$
in it.

Pick an arbitrary $x\in{{\omega}_3}\setm\sup (N_\oot\cap {{\omega}_3})$.

By transfinite recursion on ${\alpha}$, define an increasing sequence 
$\{x_{\alpha}:{\alpha}<\oot\}\subs {{\omega}_3}$ with 
$x_{\alpha}\in N_{\alpha}$ such that 
\begin{equation}
\tag{$*$}\forall n\forall a\in \br \{x_{\beta}:{\beta}<{\alpha}\};<{\omega};
g_n(a,x_{\alpha})=g_n(a,x).
\end{equation}
This can be done by $\br N_{\alpha};\oo;\subs N_{\alpha}\prec \hcal({\tau})$
and $x\in{{\omega}_3}\setm N_{\alpha}$.

Color the elements of $\{x_{\alpha}:{\alpha}<\oot\}$ with the 
$\hcal$-colors of the pair $\{x_{\alpha},x\}$:
$$
F(x_{\alpha})=\<h_n(x_{\alpha},x):n<{\omega}\>.
$$
The range of $F$ has cardinality $\le \oo^{\omega}=\oo$ by CH, so 
there is an $F$-homogeneous $D'\subs\{x_{\alpha}:{\alpha}<\oot\}$
of size $\oot$, and  it is easily seen that $D=D'\cup\{x\}$ satisfies the
 requirement of the sublemma.
\Eproof

To simplify our notation we will assume that $D$ is just 
$\oot+1$.

\begin{sublemma}
$\forall {\xi}<\oot$ $\exists {\delta}({\xi})<\oot$
$\{\fxe:{\eta}\in \oot+1\setm {\delta}({\xi})\}$ forms a $\Delta$-system.
\end{sublemma}

\Proof{sublemma}

For ${\xi}<{\eta}<{\eta}'\le\oot$ let 
$d({\xi},{\eta},{\eta}')=\dom(\fxe)\cap\dom(\fxep)$.
Then $d({\xi},{\eta},{\eta}')$ is one of the ${\omega}_1$ many countable 
subsets of $\dom(\fxe)$, so, 
by the end-homogeneity, $d({\xi},{\eta},{\eta}')$ is independent of ${\eta}'$
for ${\xi}<{\eta}<{\eta}'\le\oot$.
Denote this common value by $\hxe$.
Unfortunately we can't apply the end-homogeneity for $\hxe$, because
its range may have large cardinality.
But we can argue in the following way.
We know that ${\eta}<{\eta}'$ implies $\hxe\subs\hxep$,
because for any ${\eta}'<{\eta}''$ we have
$\dom(\fxe)\cap\dom(\fxepp)=\hxe=\dom(\fxe)\cap\dom(\fxep)$,
so $\hxe\subs \dom(\fxep)\cap \dom(\fxepp)=\hxep$.

But $\hxe$ is a countable set, so 
$$
\forall {\xi}<\oot
\ \exists \dx<\oot\ \forall {\eta},{\eta}'\ge \dx
\ \hxe=\hxep.
$$
But this means that $\{\dom(\fxe):{\eta}>\dx\}$ forms a $\Delta$-system
with kernel $h({\xi},\dx)$.
But $\fxep\rest h({\xi},\dx)$ is independent of  ${\eta}'$ by the
end homogeneity, so the functions $\fxep$ are pairwise compatible 
for $\dx\le {\eta}'\le \oot$.
\Eproof
Applying in four steps this sublemma, a $\Delta$-system argument, 
transfinite recursion and 
a $\Delta$-system argument again, we can find
a  set $D\in\br \oot+1;\oot;$ with $\oot\in D$ such that
\begin{enumerate}\arablabel
\item $\forall {\xi}\in D$ $\{\fxe:{\eta}\in D\setm {\xi}+1\}$ forms
a $\Delta$-system with kernel $\pdx$. We write
$\fxe=\pdx\cupd\rpxe$,
\item $\{\pdx:{\xi}\in D\}$ forms a $\Delta$-system with kernel $p$.
We write $\pdx=p\cupd\px$. So $\fxe=p\cupd \px\cupd \rpxe$.
\item $\forall {\eta}\in D$ $\forall {\xi},{\xi}',{\eta}'\in D\cap {\eta}$
if ${\xi}'<{\eta}'$ then 
$$
\big(\dom\pe\cup\dom\rpxe\big)\cap \dom\fxpep=\empt,
$$
\item $\{f({\xi},\oot):{\xi}\in D\}$ forms a $\Delta$-system.
\end{enumerate}
By the end-homogeneity of $D$, it follows that 
$\{\fxe:{\xi}\in D\cap {\eta}\}$ also forms a $\Delta$-system
with some kernel $\qde$ for all ${\eta}\in D$.
Write $\qde=p\cup \qe$.

Consider the increasing enumeration $\{{\delta}_{\nu}:{\nu}<\oot\}$
of $D$ and let $C=\{{\delta}_{\nu} :{\nu}$ is  limit$\}$.

\begin{Claim}{\label{cl1} If ${\xi}'<{\eta}'<{\xi}<{\eta}$ are from $D$,
then $\dom\fxpep\cap\dom\fxe=\dom(p)$}
\end{Claim}

\noindent
Indeed, $\fxe=p\cupd\px\cupd\rpxe$ and both
$\dom(\px)$ and $\dom\rpxe$ are disjoint from $\dom(\fxpep)$ by $(3)$.

\begin{Claim}
{\label{cl2}
If ${\xi}'<{\eta}'<{\eta}$ are from $C$, then
$\dom\fxpep\cap\dom\qe=\empt$
}
\end{Claim}

\noindent
Indeed, pick ${\xi}\in D\cap ({\eta}',{\eta})$, observe
$\qe\subs\fxe$, and apply claim 1.

\begin{Claim}{\label{cl3}
If ${\eta}'<{\eta}$ are from $C$, then $\dom(\qep)\cap\dom(\qe)=\empt$.
}
\end{Claim}

\noindent
Let ${\xi}'\in D\cap {\nu}'$. Then $\qep\subs \fxpep$ and apply
claim \ref{cl2}.

\begin{Claim}{\label{cl4}
If ${\xi}<{\eta}$ are from $C$, then 
$\dom(\px)\cap\dom(\qe)=\empt$.
}
\end{Claim}

\noindent
Indeed, pick ${\eta}^*<{\xi}^*\in D\cap ({\xi},{\eta})$.
Then
$\dom(\px)\cap\dom(\qe)\subs
\dom f({\xi},{\eta}^*)\cap\dom(f({\xi}^*,{\eta})\setm p))=\empt$
by claim 1.

So if you take 
$\rxe=\fxe\setm \big( \px\cup\qe\cup p\big)$ 
for ${\xi}<{\eta}\in C$, then the set $C$
and the conditions 
$\{p,\px,\qe,\rxe:{\xi},{\eta}\in C, {\xi}<{\eta}\}$ satisfy
(i)--(v). The lemma is proved.
\Eproof

\Proof{theorem \ref{tm:zoli}}

We start with a model of ZFC + GCH and  fix an arbitrarily large
regular cardinal ${\kappa}\ge {\omega}_3$. Let $\coo=\fn({\omega}_1,2,{\omega})$ and 
 consider the poset  $Q=\npr {\kappa};{\omega}_1;{\omega};\coo=
\<\fn({\kappa},\coo;{\omega}_1),\le_{\omega}\>$ where 
$f\le_{\omega}g$ if and only if  $\dom(f)\supset\dom(g)$, 
$f({\alpha})\le_{\coo} g({\alpha})$ for each
${\alpha}\in\dom(g)$ and 
$|\{{\alpha}\in\dom(g):f({\alpha})\ne g({\alpha})\}|<{\omega}$
(see \cite{FSS}).
We will show that 
the model $V^Q$ satisfies our requirements. 

It is  known that $V^Q\models 2^{\omega}={\kappa}$, $Q$ is proper
and the cardinals in $V$ and in $V^Q$ are the same
(see \cite{FSS}).
Obviously the next lemma concludes the proof.

\begin{lemma}
$V^Q\models$ ``$\forall g:\br{\omega}_3;2;\to\oo$ 
$\exists A,B\subs{\omega}_3$ $\tp(A)=\tp(B)=\oot$, $\sup A=\sup B$ and
$\exists {\alpha}\in \oo$ $\forall S\in \br A;{\omega};$
$\exists B_S\in\br B;\oot;$ $g''[S,B_S]=\{{\alpha}\}$''.
\end{lemma}

\Proof{lemma}
Assume $\gdot$ is a name of a function from 
$\br {\omega}_3;2;$ to $\oo$. 

For each ${\xi}<{\eta}<{\omega}_3$ pick a condition 
$\sxe\in Q$ and an ordinal $\axe\in {\omega}_1$
with $\sxe\force$ ``$ \gdot({\xi},{\eta})=\axe $''.
We can assume that $\dom(\sxe)\subs{\omega}_3$.

Fix an enumeration of $\coo$ in V, 
$\{s_{\nu}:{\nu}<{\omega}_1\}$, and define a bijection $F$
between $\npr {\omega}_3;{\omega}_1;{\omega};\coo;$ and $\fnn$ as follows:
$$
\begin{array}{rcl}
F(s)=h&\mbox{ if and only if  }&\dom(h)=s\\
&&\mbox{and }  s({\nu})=c_{h({\nu})}\mbox{ for all ${\nu}\in \dom (h)$}.
\end{array}
$$
Consider now the function
$f:\br {\omega}_3;2; \to P$ defined by the formula
$\fxe=F(\sxe)\times\{\<{\omega}_3,\axe\>\}$.
Formally, the range of $f$ is $P\times({\omega}_3\times \oo)=
\fn({\omega}_3+1,{\omega_1};{\omega_1})$, but this poset
is isomorphic to $P$.
Applying lemma \ref{lm:zoli} we can get
a $C\in[{\omega}_3]^{\omega_2}$ and  elements 
$\{p,\px,\qe,\rxe:{\xi},{\eta}\in C,{\xi}<{\eta}\}\subs Q$ 
and ${\alpha}\in\oo$ satisfying (i)--(v) above and $\axe={\alpha} $
for each ${\xi}< {\eta}\in C$.
Let $\adot$ and $\bdot$ be $Q$-names of subsets $C$ such that for all 
${\xi}\in C$  we have
$
\boole{{\xi}\in\adot}=F^{-1}(p\cup\px)
$
and $\boole{{\eta}\in\bdot}=F^{-1}(p\cup\qe)$.
It is clear that $F^{-1}(p)\force |\adot|=|\bdot|=\oot$.

Let $S$ be a countable subset of $A$ in $V^Q$. Since $Q$ is proper,
 there is a countable $T$ in $V$ with $S\subs T$.

Let $B^*=\{{\eta}\in B:g''[T\cap A,\{{\eta}\}]=\{{\alpha}\}\}$ and
$\bdot^*$ be a name for this set.

It is enough to show that 
$F^{-1}(p)\force |B^*|=\oot$. Assume on the contrary that 
$r\le F^{-1}(p)$, ${\rho} \in C$ and 
$r\force$ ``$B^*\subs B\cap {\rho}$'', i.e., $B^*$ is bounded in $B$.

Let $E=\dom(p)\cup\dom(r)\cup\bigcup\{\dom(p({\xi})):{\xi}\in T\}$.
Pick ${\sigma}\in C\setm {\rho}$ such that 
$\dom(q({\sigma})\cup\dom(r({\xi},{\sigma})\cap E=\empt$ for
each ${\xi}\in T$. Since $r\force |S|={\omega}$, there is ${\xi}\in T$
such that $r$ and $p\cup p({\xi})$ is compatible.

Let $r^*=r\land (p\cup p({\xi})\land (q({\sigma})\cup r({\xi},{\sigma})$.
Then $r^*$ forces a contradiction.
\Eproof
This completes the proof of the theorem.
\Eproof

\section{Models without large  first countable  $P_2$-spaces}
\label{sc:cohen}

Hajnal and Juh{\'a}sz \cite[Problem 10]{HJ} asked if it is consistent
to assume that $2^{\omega}\ge {\omega_2} $ and every first countable
(or compact) space with property $P_2$ has cardinality 
$\le {\omega_1}$. We  give an affirmative answer here.
We will argue in the following way. 
First  we quote the  definition of principle $\cax({\kappa})$ from
\cite{CAX}, then we show that $\cax(\oot)$ implies that every 
 first countable $P_2$ space has cardinality $\le\oo$, finally 
we will cite a theorem from \cite{CAX} saying that $\cax(\oot)$ is 
consistent with any cardinal arithmetic. 

\begin{definition}\label{ax:c}(See \cite{CAX}.)
Let ${\kappa}$ be an infinite cardinal.
We say that principle $\cax({\kappa})$ holds if and only if  for each family $\{A({\xi},n),B({\xi},n):
{\xi}\in {\kappa},n\in {\omega}\}\subs \br {\omega};{\omega};$ 
either $($i$)$ or $($ii$)$ below holds:
\begin{enumerate}\rlabel
\item 
$\exists C\in \br {\kappa};{\kappa};$ $\forall n,m\in {\omega}$
$\forall {\xi}\ne {\zeta}\in C$ $A({\xi},n)\cap B({\zeta},m)\ne\empt$,
\item  
$\exists D,E\in \br {\kappa};{\kappa};$ $\exists n,m\in {\omega}$
$\forall {\xi}\in D$ $\forall {\zeta}\in E$ 
$A({\xi},n)\cap B({\zeta},m)=\empt$.
\end{enumerate}
\end{definition}

\begin{theorem}\label{tm:p2}
If $\cax({\kappa})$ holds then each
 first countable, separable Hausdorff space $X$ of size ${\kappa}$
contains two disjoint open sets $U$ and $V$ of cardinality ${\kappa}$.
\end{theorem}

\proof
Let $S$ be a countable dense subset of $X$. For each $x\in X$
fix a neighborhood base $\{U(x,n):n\in {\omega}\}$ of $x$ in $X$.
Take $A(x,n)=B(x,n)=U(x,n)\cap S$ and apply 
$\cax({\kappa})$. Since $X$ is $T_2$, there is no 
$C\in\br X;{\kappa};$ satisfying \ref{ax:c}(i). So there are 
$D,E\in \br X;{\kappa};$ and $n,m\in {\omega}$ such that 
$U(x,n)\cap U(y,m)\cap S=\empt$ whenever $x\in D$ and $y\in E$.
But $S$ is dense in $X$, therefore $U=\bigcup\{U(x,n):x\in D\}$ and
$V=\bigcup\{U(y,m):y\in E\}$ are disjoint and of cardinality ${\kappa}$.
\eproof

It was proved in \cite{CAX} that  starting from a model of CH, after adding ${\lambda}$-many Cohen reals by the poset $P=\fn({\lambda},2,{\omega})$, 
we have that $\cax(\oot)$ holds in $V^P$.
Since $P_2$ spaces are separable as it was observed in the proof of
\cite[theorem 1]{HJ}, theorem \ref{tm:p2} yields the following corollary.

\begin{corollary}\label{tm:no-p2}
If ZF is consistent then so is ZFC + {\em ``$\cont$ is as large as you wish''} 
+ {\em ``every first countable $P_2$ space has cardinality 
$\le{\omega_1} $.''}
\end{corollary}


\begin{thebibliography}{99}
\bibitem{BS} J. E. Baumgartner, S. Shelah, {\em  Remarks on
superatomic Boolean algebras}, Ann. Pure Appl. Logic, {\bf 33}
(l987), no. 2, 1-9-129.
\bibitem{FSS} S. Fuchino, S. Shelah, L. Soukup, 
{\em Beating with sticks and clubs}, 1994, preprint.
\bibitem{H} A. Hajnal, {\em Proof of a conjecture of S. Ruziewicz},
Fund. Math. {\bf 50} (1961), 123-128.
\bibitem{HJ} A. Hajnal, I. Juh{\'a}sz,
{\em Intersection properties of open sets}, Topology and its
Application,
{\bf 19} (l985), 201--209.
\bibitem{He}  S. H. Hechler, {\em Short complete nested sequences in
$\beta N - N$ and small maximal almost disjoint families},
Gen. Top. Appl. {\bf 2} (1972) 139-149.
\bibitem{J} I. Juh{\'a}sz,  
{\em Cardinal Functions - Ten Years Later}, 
Math. Center Tracts 123, Amsterdam, 1980.
 \bibitem{CAX} I.  Juh{\'a}sz, S. Soukup, Z. Szentmikl{\'o}ssy, 
{\em Combinatorial principles from  adding Cohen reals}, to appear 
in  Proc. of Logic Colloquium 95.
\bibitem{vM} J. Van Mill {\em An introduction to ${\beta}N$}, 
in Handbook of Set-Theoretic Topology, 
Ed. K. Kunen J.E Vaughan, Elsevier Publ. B. V. 1984.
\bibitem{Sh} S. Shelah,  private communication, 1993.
\bibitem{Sh2} S. Shelah, {\em Models with second order properties III.
Omitting types in ${\lambda}^+$ for $L(Q)$},
Proc. of a Berlin Workshop July 1977, 
Archive f. Math. Logik {\bf 12} (l980) 1-11.
\bibitem{V} D. Velleman {\em ${\omega}$-morasses and a weak form
of Martin's axiom provable in ZFC}, 
Trans. Amer. Math. Soc. {\bf 285} (1984), no 2. 617--627
\end{thebibliography}
\end{document}